\documentclass{amsart}
\usepackage{upref}

\theoremstyle{definition}

\theoremstyle{remark}

\begin{document}

\title{Real and Complex Integrals on Spheres and Balls}

\author{H. Turgay Kaptano\u{g}lu}
\address{Department of Mathematics, B$\dot{\hbox{\i}}$lkent University,
Ankara 06800, Turkey}
\email{kaptan@fen.bilkent.edu.tr}
\urladdr{http://www.fen.bilkent.edu.tr/$\sim$kaptan/}

\date{13 February 2014}


\begin{abstract}
We evaluate integrals of certain polynomials over spheres and balls in real
or complex spaces.
We also promote the use of the Pochhammer symbol which gives the values of
our integrals in compact forms.
\end{abstract}

\maketitle

\section{Introduction}

The integrals we are interested in arise in the computations of certain norms
in weighted Bergman and Besov spaces of harmonic or holomorphic functions on
the unit real or complex ball.
They are
$$\int_{\mathbb B}|u^\alpha|^p\,(1-|u|^2)^q\,dV(u)\qquad\text{and}\qquad
\int_{\mathbb B}|\langle u,c\rangle|^p\,(1-|u|^2)^q\,dV(u)$$
for $p\ge0$ and $q>-1$ or their variants.
We compute their exact values in terms of the gamma function.
For $q\le-1$, the integrals are clearly divergent.

The values of these integrals are used quite often in practical computations
such as in obtaining maximal growth rates for norms or function values in
the spaces mentioned, or when the norms of the same function in different
spaces are compared, or in constructing a function in a given space that does
not lie in another space or with prescribed growth properties.

Let's start by explaining the notation.
We work in $\mathbb R^n$ or $\mathbb C^N$, and identify even-dimensional
spaces by always equating $2N=n$ with $n\ge2$ and $N\ge1$.
Arbitrary points are $x=(x_1,\ldots,x_n)\in\mathbb R^n$ and
$z=(z_1,\ldots,z_N)\in\mathbb C^N$.
We show the inner product and its norm on either space is by
$\langle\cdot,\cdot\rangle$ and $|\cdot|$.
So $\langle x,y\rangle=x_1y_1+\cdots+x_ny_n$,
$\langle z,w\rangle=z_1\overline w_1+\cdots+z_N\overline w_N$,
$|x|^2=x_1^2+\cdots+x_n^2$, and $|z|^2=|z_1|^2+\cdots+|z_N|^2$.
We use $\mathbb B$ and $\mathbb S$ to denote the unit ball and the unit
sphere with respect to these norms; so $\mathbb B=\{\,|\cdot|<1\}$ and
$\mathbb S=\{\,|\cdot|=1\}$.
We also write $x=r\xi$ or $z=r\zeta$ with $r=|x|$ or $r=|z|$, and
$\xi,\zeta\in\mathbb S$; that is, $(x_1,\ldots,x_n)=r\,(\xi_1,\ldots,\xi_n)$,
and $(z_1,\ldots,z_n)=r\,(\zeta_1,\ldots,\zeta_n)$,

We use multi-index notation in which $\alpha$ is an $n$- or $N$-tuple of
nonnegative integers.
So in the real case, $\alpha=(\alpha_1,\ldots,\alpha_n)$,
$|\alpha|=\alpha_1+\cdots+\alpha_n$, $\alpha!=\alpha_1!\cdots\alpha_n!$,
$0^0=1$, and $x^\alpha=x_1^{\alpha_1}\cdots x_n^{\alpha_n}$.
Then
\begin{align*}
|x^\alpha|^p
&=|x_1^{\alpha_1}\cdots x_n^{\alpha_n}|^p
=|x_1|^{\alpha_1p}\cdots|x_n|^{\alpha_np}\\
&=r^{\alpha_1p}\,|\xi_1|^{\alpha_1p}\cdots r^{\alpha_np}\,|\xi_n|^{\alpha_np}
=r^{|\alpha|p}\,|\xi^\alpha|^p.
\end{align*}
In the complex case, $\alpha=(\alpha_1,\ldots,\alpha_N)$ and
$|z^\alpha|^p=|z_1^{\alpha_1}\cdots z_N^{\alpha_N}|^p
=r^{|\alpha|p}\,|\zeta^\alpha|^p$.

We denote the usual multidimensional Lebesgue measure by $V$ and the
one-lower-dimensional surface measure it induces on $\mathbb S$ by $S$.
Their traditional forms are such that the unit cube and the unit square have
measure $1$.
We also use a normalization of these measures by letting $\nu$ and $\sigma$
denote the volume and surface measures with $\nu(\mathbb B)=1$ and
$\sigma(\mathbb S)=1$.

It is well-known that (see \cite{F})
$$V(\mathbb B)=\frac{\pi^{n/2}}{\Gamma(1+n/2)}\qquad\text{and}\qquad
S(\mathbb S)=\frac{n\,\pi^{n/2}}{\Gamma(1+n/2)},$$
where
$$\Gamma(t)=\int_0^\infty s^{t-1}\,e^{-s}\,ds
=2\int_0^\infty r^{2t-1}\,e^{-r^2}\,dr\qquad(t>0)$$
is the \textit{gamma function}.
So
$$d\nu=\frac{\Gamma(1+n/2)}{\pi^{n/2}}\,dV\qquad\text{and}\qquad
d\sigma=\frac{\Gamma(n/2)}{2\,\pi^{n/2}}\,dS.$$
Taking $n=2N$, we obtain the values in $\mathbb C^N$ as
$$d\nu=\frac{N!}{\pi^N}\,dV\qquad\text{and}\qquad
d\sigma=\frac{(N-1)!}{2\,\pi^N}\,dS.$$

There is also the related beta integral
$$2\int_0^1r^{2a-1}\,(1-r^2)^{b-1}\,dr=\int_0^1s^{a-1}\,(1-s)^{b-1}\,ds
=\mathrm B(a,b)=\frac{\Gamma(a)\,\Gamma(b)}{\Gamma(a+b)}$$
for $a,b>0$.
Let's further recall that $\Gamma(t+1)=t\,\Gamma(t)$, so $\Gamma(k)=(k-1)!$
for positive integer $k$, $\Gamma(2)=\Gamma(1)=1$, and
$\Gamma(1/2)=\sqrt\pi$.

The polar coordinates formula
$$\int_{\mathbb B}f(x)\,dV(x)
=\int_0^1r^{n-1}\int_{\mathbb S}f(r\xi)\,dS(\xi)\,dr$$
takes the form
$$\int_{\mathbb B}f(x)\,d\nu(x)
=n\int_0^1r^{n-1}\int_{\mathbb S}f(r\xi)\,d\sigma(\xi)\,dr$$
with the normalization.
If the integration is on all of $\mathbb R^n$, then the upper limit of
integration with respect to $r$ is $\infty$.
Simply taking $n=2N$ yields their forms in $\mathbb C^N$; see
\cite[1.4.3]{R}.

The \textit{Pochhammer symbol} $(a)_b$ is defined by
$$(a)_b=\frac{\Gamma(a+b)}{\Gamma(a)}$$
whenever the gamma functions make sense.
This is a shifted factorial, because $(a)_k=a(a+1)\cdots(a+k-1)$ for positive
integer $k$ and for all $a\in\mathbb C$.
In particular, $(1)_k=k!$ and $(a)_0=1$.
This function includes the gamma function as a special case and thus is more
useful.
It has a slower growth than the factorial; see \cite{TE}.
Stirling formula gives
\begin{equation}\label{growth}
\frac{\Gamma(c+a)}{\Gamma(c+b)}\sim c^{a-b},\qquad
\frac{(a)_c}{(b)_c}\sim c^{a-b},\qquad
\frac{(c)_a}{(c)_b}\sim c^{a-b},\qquad(c\to\infty),
\end{equation}
where $u\sim v$ means that $|u/v|$ is bounded above and below by two positive
constants.

Finally, we use $:=$ to mean that the right hand side defines the left hand
side.

\section{Real Integrals}

As previously utilized in \cite{F} and \cite[1.4.9]{R}, we start with
$$J_1:=\int_{\mathbb R^n}|x^\alpha|^p\,e^{-|x|^2}\,dV(x)$$
and evaluate it in two ways; first by writing it as a product of
one-dimensional integrals, and second by using the polar coordinates formula.
First,
$$J_1=\prod_{j=1}^n\int_{-\infty}^\infty|x_j|^{\alpha_jp}\,e^{-x_j^2}\,dx_j
=\prod_{j=1}^n2\int_0^\infty x_j^{\alpha_jp}\,e^{-x_j^2}\,dx_j,$$
where each integral is a gamma function.
So 
$$J_1=\int_{\mathbb R^n}|x^\alpha|^p\,e^{-|x|^2}\,dV(x)
=\prod_{j=1}^n\Gamma\Bigl(\frac{1+\alpha_jp}2\Bigr).$$
Second,
\begin{align*}
J_1&=\int_0^\infty r^{n-1}\,r^{|\alpha|p}\,e^{-r^2}
\int_{\mathbb S}|\xi^\alpha|^p\,dS(\xi)\,dr\\
&=J_2\int_0^\infty r^{n-1+|\alpha|p}\,e^{-r^2}\,dr
=\frac12\,J_2\,\Gamma\Bigl(\frac{n+|\alpha|p}2\Bigr),
\end{align*}
where
$$J_2:=\int_{\mathbb S}|\xi^\alpha|^p\,dS(\xi)$$
Comparing the two forms of $J_1$, we obtain the value
$$J_2=\int_{\mathbb S}|\xi^\alpha|^p\,dS(\xi)
=\frac{2\prod_{j=1}^n\Gamma\bigl((1+\alpha_jp)/2\bigr)}
{\Gamma\bigl((n+|\alpha|p)/2\bigr)}.$$

Continuing, directly using the polar coordinates formula and then the beta
integral, we obtain
\begin{align*}
\int_{\mathbb B}|x^\alpha|^p\,(1-|x|^2)^q\,dV(x)
&=\int_0^1r^{n-1}\,r^{|\alpha|p}\,(1-r^2)^q
\int_{\mathbb S}|\xi^\alpha|^p\,dS(\xi)\,dr\\
&=J_2\int_0^1r^{n-1+|\alpha|p}\,(1-r^2)^q\,dr\\
&=\frac{2\prod_{j=1}^n\Gamma\bigl((1+\alpha_jp)/2\bigr)}
{\Gamma\bigl((n+|\alpha|p)/2\bigr)}\,\frac12
\,\frac{\Gamma\bigl((n+|\alpha|p)/2\bigr)\,\Gamma(1+q)}
{\Gamma\bigl((n+|\alpha|p)/2+1+q\bigr)}.
\end{align*}
Hence
$$J_3:=\int_{\mathbb B}|x^\alpha|^p\,(1-|x|^2)^q\,dV(x)
=\frac{\Gamma(1+q)\,\prod_{j=1}^n\Gamma\bigl((1+\alpha_jp)/2\bigr)}
{\Gamma\bigl(1+q+(n+|\alpha|p)/2\bigr)}.$$
Letting $q=0$ gives
$$J_4:=\int_{\mathbb B}|x^\alpha|^p\,dV(x)
=\frac{\prod_{j=1}^n\Gamma\bigl((1+\alpha_jp)/2\bigr)}
{\Gamma\bigl(1+(n+|\alpha|p)/2\bigr)},$$
which is the volume counterpart of $J_2$.

We are also interested in the versions of $J_1$--$J_4$ in which the
integration is with respect to $\sigma$ and $\nu$.
First of all,
$$\frac{\prod_{j=1}^n\Gamma(1/2+\alpha_jp/2)}{\pi^{n/2}}
=\frac{\prod_{j=1}^n\Gamma(1/2)\,(1/2)_{\alpha_jp/2}}{\Gamma(1/2)^n}
=\prod_{j=1}^n\Bigl(\frac12\Bigr)_{\alpha_jp/2}.$$
Then we have
$$\int_{\mathbb R^n}|x^\alpha|^p\,e^{-|x|^2}\,d\nu(x)
=\frac{\Gamma(1+n/2)}{\pi^{n/2}}\,J_1
=\Gamma\Bigl(1+\frac n2\Bigr)
\,\frac{\prod_{j=1}^n\Gamma(1/2+\alpha_jp/2)}{\Gamma(1/2)^n}$$
and
$$J_1':=\int_{\mathbb R^n}|x^\alpha|^p\,e^{-|x|^2}\,d\nu(x)
=(1)_{n/2}\,\prod_{j=1}^n\Bigl(\frac12\Bigr)_{\alpha_jp/2}.$$
Next,
$$\int_{\mathbb S}|\xi^\alpha|^p\,d\sigma(\xi)
=\frac{\Gamma(n/2)}{2\,\pi^{n/2}}\,J_2
=\frac{\Gamma(n/2)}{2\,\Gamma(1/2)^n}
\,\frac{2\prod_{j=1}^n\Gamma(1/2+\alpha_jp/2)}{\Gamma(n/2+|\alpha|p/2)}$$
and
$$J_2':=\int_{\mathbb S}|\xi^\alpha|^p\,d\sigma(\xi)
=\frac{\prod_{j=1}^n(1/2)_{\alpha_jp/2}}{(n/2)_{|\alpha|p/2}}.$$
Finally,
\begin{align*}
\int_{\mathbb B}|x^\alpha|^p\,(1-|x|^2)^q\,d\nu(x)
&=\frac{\Gamma(1+n/2)}{\pi^{n/2}}\,J_3\\
&=\frac{\Gamma(1+n/2)}{\Gamma(1/2)^n}
\,\frac{\Gamma(1+q)\,\prod_{j=1}^n\Gamma(1/2+\alpha_jp/2)}
{\Gamma(1+q+n/2+|\alpha|p/2)},
\end{align*}
$$J_3':=\int_{\mathbb B}|x^\alpha|^p\,(1-|x|^2)^q\,d\nu(x)
=\frac{(1)_q\,\prod_{j=1}^n(1/2)_{\alpha_jp/2}}{(1+n/2)_{q+|\alpha|p/2}},$$
and
$$J_4':=\int_{\mathbb B}|x^\alpha|^p\,d\nu(x)
=\frac{\prod_{j=1}^n(1/2)_{\alpha_jp/2}}{(1+n/2)_{|\alpha|p/2}}.$$

It is interesting that the values of integrals using the usual Lebesgue
measures have more natural and simpler appearance with the gamma function
while the values of integrals using the normalized Lebesgue measures with the
Pochhammer symbol.

The integrals assume especially simple forms when $p=2m$ is a positive even
integer and $q=k$ is a positive integer.
In such a case $|x^\alpha|^{2m}=x^{2m\alpha}$, where
$2m\alpha=(2m\alpha_1,\ldots,2m\alpha_n)$, and $\Gamma(1+k)=(1)_k=k!$.
Thus
\begin{align*}
J_1''&:=\int_{\mathbb R^n}x^{2m\alpha}\,e^{-|x|^2}\,dV(x)
=\prod_{j=1}^n\Gamma\Bigl(\frac12+m\alpha_j\Bigr),\\
J_2''&:=\int_{\mathbb S}\xi^{2m\alpha}\,dS(\xi)
=\frac{2\prod_{j=1}^n\Gamma(1/2+m\alpha_j)}{\Gamma(n/2+m|\alpha|)},\\
J_3''&:=\int_{\mathbb B}x^{2m\alpha}\,(1-|x|^2)^k\,dV(x)
=\frac{k!\,\prod_{j=1}^n\Gamma(1/2+m\alpha_j)}
{\Gamma(1+k+n/2+m|\alpha|)},\\
J_4''&:=\int_{\mathbb B}x^{2m\alpha}\,dV(x)
=\frac{\prod_{j=1}^n\Gamma(1/2+m\alpha_j)}{\Gamma(1+n/2+m|\alpha|)},
\end{align*}
and
\begin{align*}
J_1'''&:=\int_{\mathbb R^n}x^{2m\alpha}\,e^{-|x|^2}\,d\nu(x)
=(1)_{n/2}\,\prod_{j=1}^n\Bigl(\frac12\Bigr)_{m\alpha_j},\\
J_2'''&:=\int_{\mathbb S}\xi^{2m\alpha}\,d\sigma(\xi)
=\frac{\prod_{j=1}^n(1/2)_{m\alpha_j}}{(n/2)_{m|\alpha|}},\\
J_3'''&:=\int_{\mathbb B}x^{2m\alpha}\,(1-|x|^2)^k\,d\nu(x)
=\frac{k!\,\prod_{j=1}^n(1/2)_{m\alpha_j}}{(1+n/2)_{k+m|\alpha|}},\\
J_4'''&:=\int_{\mathbb B}x^{2m\alpha}\,d\nu(x)
=\frac{\prod_{j=1}^n(1/2)_{m\alpha_j}}{(1+n/2)_{m|\alpha|}}.
\end{align*}

Further simplifications are possible in the case of integer parameters by
reducing all gamma functions to $\Gamma(1)$ or $\Gamma(1/2)$ and by treating
even $n$ and odd $n$ separately, and then rewriting everything in terms of
factorials and $\sqrt{\pi}$, but we prefer to stop here.
We have the occasion to do so in complex integrals.

Note that the integrals containing $\xi^\alpha$ or $x^\alpha$ without
absolute values and with an $\alpha$ not consisting of all even entries is
$0$ by symmetry; see \cite{F}.

\section{Real Integrals Depending on a Parameter}

Our next set of integrals contain inner products with a second variable
$y$, but the integration is with respect to the same variable as before.
Consequently our formulas are expressed in terms of $y$ yielding the growth
rate of the integral as a function of $|y|$.

Lebesgue measures on the sphere or on the ball, normalized or not,  are
rotationally invariant.
Thus in our computations, there is no loss of generality to fix a simple
$y=y_1(1,0,\ldots,0)$ so that $|y|^2=y_1^2$, $\langle\xi,y\rangle=\xi_1y_1$,
and $|\langle\xi,y\rangle|=|\xi_1|\,|y|$.

We start with
$$\int_{\mathbb S}|\langle\xi,y\rangle|^p\,dS(\xi)
=|y|^p\int_{\mathbb S}|\xi_1|^p\,dS(\xi)=|y|^p\,J_2,$$
where $J_2$ is in the special case with $\alpha=(1,0,\ldots,0)$ and
$|\alpha|=1$.
Recalling that $\Gamma(1/2)=\sqrt{\pi}$,
$$J_6:=\int_{\mathbb S}|\langle\xi,y\rangle|^p\,dS(\xi)
=\frac{2\,\pi^{(n-1)/2}\,\Gamma\bigl((1+p)/2\bigr)}
{\Gamma\bigl((n+p)/2\bigr)}\,|y|^p.$$

Next,
\begin{align*}
\int_{\mathbb R^n}|\langle x,y\rangle|^p\,e^{-|x|^2}\,dV(x)
&=\int_0^\infty r^{n-1}\,r^p\,e^{-r^2}\int_{\mathbb S}|\langle\xi,y\rangle|^p
\,dS(\xi)\,dr\\
&=J_6\int_0^\infty r^{n-1}\,r^p\,e^{-r^2}\,dr\\
&=\frac{2\,\pi^{(n-1)/2}\,\Gamma\bigl((1+p)/2\bigr)}
{\Gamma\bigl((n+p)/2\bigr)}\,|y|^p\,\frac12\,\Gamma\Bigl(\frac{n+p}2\Bigr),
\end{align*}
and hence
$$J_5:=\int_{\mathbb R^n}|\langle x,y\rangle|^p\,e^{-|x|^2}\,dV(x)
=\pi^{(n-1)/2}\,\Gamma\Bigl(\frac{1+p}2\Bigr)\,|y|^p.$$

Further,
\begin{align*}
\int_{\mathbb B}|\langle x,y\rangle|^p\,(1-|x|^2)^q\,dV(x)
&=\int_0^1r^{n-1}\,r^p\,(1-r^2)^q\int_{\mathbb S}|\langle\xi,y\rangle|^p
\,dS(\xi)\,dr\\
&=J_6\int_0^1r^{n-1+p}\,(1-r^2)^q\,dr\\
&=\frac{2\,\pi^{(n-1)/2}\,\Gamma\bigl((1+p)/2\bigr)}
{\Gamma\bigl((n+p)/2\bigr)}\,|y|^p
\,\frac12\,\frac{\Gamma\bigl((n+p)/2\bigr)\,\Gamma(1+q)}
{\Gamma\bigl((n+p)/2+1+q\bigr)}
\end{align*}
with the help of the beta integral.
Thus
$$J_7:=\int_{\mathbb B}|\langle x,y\rangle|^p\,(1-|x|^2)^q\,dV(x)
=\frac{\pi^{(n-1)/2}\,\Gamma(1+q)\,\Gamma\bigl((1+p)/2\bigr)}
{\Gamma\bigl(1+q+(n+p)/2\bigr)}\,|y|^p,$$
and when $q=0$,
$$J_8:=\int_{\mathbb B}|\langle x,y\rangle|^p\,dV(x)
=\frac{\pi^{(n-1)/2}\,\Gamma\bigl((1+p)/2\bigr)}
{\Gamma\bigl(1+(n+p)/2\bigr)}\,|y|^p.$$

The versions of $J_6$--$J_8$ with $\sigma$ and $\nu$ are computed just
as in $J_1'$--$J_4'$ and are
\begin{align*}
J_5'&:=\int_{\mathbb R^n}|\langle x,y\rangle|^p\,e^{-|x|^2}\,d\nu(x)
=\Gamma\Bigl(1+\frac n2\Bigr)\,\Bigl(\frac12\Bigr)_{\!p/2}\,|y|^p,\\
J_6'&:=\int_{\mathbb S}|\langle\xi,y\rangle|^p\,d\sigma(\xi)
=\frac{(1/2)_{p/2}}{(n/2)_{p/2}}\,|y|^p,\\
J_7'&:=\int_{\mathbb B}|\langle x,y\rangle|^p\,(1-|x|^2)^q\,d\nu(x)
=\frac{(1)_q\,(1/2)_{p/2}}{(1+n/2)_{q+p/2}}\,|y|^p,\\
J_8'&:=\int_{\mathbb B}|\langle x,y\rangle|^p\,d\nu(x)
=\frac{(1/2)_{p/2}}{(1+n/2)_{p/2}}\,|y|^p.
\end{align*}
Note that the values of these latter integrals with respect to normalized
measures have much simpler forms.

Again let's write down the same integrals when the parameters are integers as
$p=2m$ and $q=k$.
We have
\begin{align*}
J_5''&:=\int_{\mathbb R^n}\langle x,y\rangle^{2m}\,e^{-|x|^2}\,dV(x)
=\pi^{(n-1)/2}\,\Gamma\Bigl(\frac12+m\Bigr)\,|y|^{2m},\\
J_6''&:=\int_{\mathbb S}\langle\xi,y\rangle^{2m}\,dS(\xi)
=\frac{2\,\pi^{(n-1)/2}\,\Gamma(1/2+m)}{\Gamma(n/2+m)}\,|y|^{2m},\\
J_7''&:=\int_{\mathbb B}\langle x,y\rangle^{2m}\,(1-|x|^2)^k\,dV(x)
=\frac{\pi^{(n-1)/2}\,k!\,\Gamma(1/2+m)}{\Gamma(1+k+n/2+m)}\,|y|^{2m},\\
J_8''&:=\int_{\mathbb B}\langle x,y\rangle^{2m}\,dV(x)
=\frac{\pi^{(n-1)/2}\,\Gamma(1/2+m)}{\Gamma(n/2+m)}\,|y|^{2m},
\end{align*}
and
\begin{align*}
J_5'''&:=\int_{\mathbb R^n}\langle x,y\rangle^{2m}\,e^{-|x|^2}\,d\nu(x)
=\Gamma\Bigl(1+\frac n2\Bigr)\,\Bigl(\frac12\Bigr)_{\!m}\,|y|^{2m},\\
J_6'''&:=\int_{\mathbb S}\langle\xi,y\rangle^{2m}\,d\sigma(\xi)
=\frac{(1/2)_m}{(n/2)_m}\,|y|^{2m},\\
J_7'''&:=\int_{\mathbb B}\langle x,y\rangle^{2m}\,(1-|x|^2)^k\,d\nu(x)
=\frac{k!\,(1/2)_m}{(1+n/2)_{k+m}}\,|y|^{2m},\\
J_8'''&:=\int_{\mathbb B}\langle x,y\rangle^{2m}\,d\nu(x)
=\frac{(1/2)_m}{(1+n/2)_m}\,|y|^{2m}.
\end{align*}

\section{Complex Integrals}

It is not enough just to substitute $n=2N$ to evaluate the integrals in
$\mathbb C^N$, because $z^\alpha$ is not a special case of $x^\alpha$.
This difference is more visible in the computation of $K_1$ which replaces
$J_1$ above.
Here even the one-dimensional integrals require the classical polar
coordinates formula.
So
\begin{align*}
K_1&:=\int_{\mathbb C^N}|z^\alpha|^p\,e^{-|z|^2}\,dV(z)
=\prod_{j=1}^N\int_{\mathbb C}|z_j|^{\alpha_jp}\,e^{-|z_j|^2}\,dV(z_j)\\
&=\prod_{j=1}^N\int_0^\infty r_j\,r_j^{\alpha_jp}\,e^{-r_j^2}
\int_{\mathbb S}1\,dS(\zeta_j)\,dr_j
=\pi^N\prod_{j=1}^N2\int_0^\infty r_j^{1+\alpha_jp}\,e^{-r_j^2}\,dr_j\\
&=\pi^N\prod_{j=1}^N\Gamma\Bigl(1+\frac{\alpha_jp}2\Bigr),
\end{align*}
where $dV(z_j)$ indicates an area integral and $\mathbb S$ is the unit
circle, and
\begin{align*}
K_1&=\int_0^\infty r^{2N-1}\,r^{|\alpha|p}\,e^{-r^2}
\int_{\mathbb S}|\zeta^\alpha|^p\,dS(\zeta)\,dr\\
&=K_2\int_0^\infty r^{2N-1+|\alpha|p}\,e^{-r^2}\,dr
=\frac12\,K_2\,\Gamma\Bigl(N+\frac{|\alpha|p}2\Bigr).
\end{align*}
Thus
$$K_2:=\int_{\mathbb S}|\zeta^\alpha|^p\,dS(\zeta)
=\frac{2\pi^N\prod_{j=1}^N\Gamma(1+\alpha_jp/2)}{\Gamma(N+|\alpha|p/2)}.$$
Further
\begin{align*}
K_3&:=\int_{\mathbb B}|z^\alpha|^p\,(1-|z|^2)^q\,dV(z)
=\int_0^1r^{2N-1}\,r^{|\alpha|p}\,(1-r^2)^q
\int_{\mathbb S}|\zeta^\alpha|^p\,dS(\zeta)\,dr\\
&=K_2\int_0^1r^{2N-1+|\alpha|p}\,(1-r^2)^q\,dr\\
&=\frac{2\pi^N\prod_{j=1}^N\Gamma(1+\alpha_jp/2)}
{\Gamma(N+|\alpha|p/2)}\,\frac12
\,\frac{\Gamma(N+|\alpha|p/2)\,\Gamma(1+q)}{\Gamma(N+|\alpha|p/2+1+q)}\\
&=\frac{\pi^N\,\Gamma(1+q)\,\prod_{j=1}^N\Gamma(1+\alpha_jp/2)}
{\Gamma(1+q+N+|\alpha|p/2)}
\end{align*}
and
$$K_4:=\int_{\mathbb B}|z^\alpha|^p\,dV(z)
=\frac{\pi^N\,\prod_{j=1}^N\Gamma(1+\alpha_jp/2)}
{\Gamma(1+N+|\alpha|p/2)}
$$

The integrals with respect to the normalized measures are obtained in exactly
the same way as in the real case, only simpler.
We have
$$
K_1':=\int_{\mathbb C^N}|z^\alpha|^p\,e^{-|z|^2}\,d\nu(z)
=\frac{N!}{\pi^N}\,K_1
=N!\prod_{j=1}^N\Gamma\Bigl(1+\frac{\alpha_jp}2\Bigr)
=N!\prod_{j=1}^n(1)_{\alpha_jp/2},$$
\begin{align*}
K_2'&:=\int_{\mathbb S}|\zeta^\alpha|^p\,d\sigma(\zeta)
=\frac{(N-1)!}{2\pi^N}\,K_2\\
&=\Gamma(N)\,\frac{\prod_{j=1}^N\Gamma(1+\alpha_jp/2)}{\Gamma(N+|\alpha|p/2)}
=\frac{\prod_{j=1}^N(1)_{\alpha_jp/2}}{(N)_{|\alpha|p/2}},
\end{align*}
\begin{align*}
K_3'&:=\int_{\mathbb B}|z^\alpha|^p\,(1-|z|^2)^q\,d\nu(z)
=\frac{N!}{\pi^N}\,K_3\\
&=\Gamma(1+N)\,\frac{\Gamma(1+q)\,\prod_{j=1}^N\Gamma(1+\alpha_jp/2)}
{\Gamma(1+q+N+|\alpha|p/2)}
=\frac{(1)_q\,\prod_{j=1}^N(1)_{\alpha_jp/2}}{(1+N)_{q+|\alpha|p/2}},
\end{align*}
and
$$K_4':=\int_{\mathbb B}|z^\alpha|^p\,d\nu(z)
=\frac{\prod_{j=1}^N(1)_{\alpha_jp/2}}{(1+N)_{|\alpha|p/2}}.$$

When $p=2m$ is a positive even integer and $q=k$ is a positive integer, the
integrals assume even simpler forms than the real integrals since the $1/2$
in the argument of the gamma function is now replaced by $1$.
This further allows us to use factorials instead.
The underlying reason of course is the even-dimensionality of the underlying
space.
But now we cannot remove the moduli on $z^\alpha$, because otherwise the
integrals are all zero; see \cite[Proposition 1.4.8]{R}.
Recalling the multi-index notation
$(m\alpha)!=(m\alpha_1)!\cdots(m\alpha_N)!$, we have
\begin{align*}
K_1''&:=\int_{\mathbb C^N}|z^{2m\alpha}|\,e^{-|z|^2}\,dV(z)
=\pi^N\,(m\alpha)!,\\
K_2''&:=\int_{\mathbb S}|\zeta^{2m\alpha}|\,dS(\zeta)
=\frac{2\pi^N\,(m\alpha)!}{(N-1+m|\alpha|)!},\\
K_3''&:=\int_{\mathbb B}|z^{2m\alpha}|\,(1-|z|^2)^k\,dV(z)
=\frac{\pi^N\,k!\,(m\alpha)!}{(N+k+m|\alpha|)!},\\
K_4''&:=\int_{\mathbb B}|z^{2m\alpha}|\,dV(z)
=\frac{\pi^N\,(m\alpha)!}{(N+m|\alpha|)!},
\end{align*}
and
\begin{align*}
K_1'''&:=\int_{\mathbb C^N}|z^{2m\alpha}|\,e^{-|z|^2}\,d\nu(z)
=N!\,(m\alpha)!,\\
K_2'''&:=\int_{\mathbb S}|\zeta^{2m\alpha}|\,d\sigma(\zeta)
=\frac{(N-1)!\,(m\alpha)!}{(N-1+m|\alpha|)!},\\
K_3'''&:=\int_{\mathbb B}|z^{2m\alpha}|\,(1-|z|^2)^k\,d\nu(z)
=\frac{N!\,k!\,(m\alpha)!}{(N+k+m|\alpha|)!},\\
K_4'''&:=\int_{\mathbb B}|z^{2m\alpha}|\,d\nu(z)
=\frac{N!\,(m\alpha)!}{(N+m|\alpha|)!}.
\end{align*}

To sum up, in passing from the real $J_1'$--$J_4'$ to the complex
$K_1'$--$K_4'$, essentially we only replace $1/2$ by $1$ with the appropriate
modification for the dimension.
In passing from $J_1$--$J_4$ to $K_1$--$K_4$, additionally coefficients of
$\pi^N$ or $2\pi^N$ have to be taken care of.

\section{Complex Integrals Depending on a Parameter}

After having already shown the details a few times, we are now content with
simply stating the results.
\begin{align*}
K_5&:=\int_{\mathbb C^N}|\langle z,w\rangle|^p\,e^{-|z|^2}\,dV(z)
=\pi^N\,\Gamma\Bigl(1+\frac p2\Bigr)\,|w|^p,\\
K_6&:=\int_{\mathbb S}|\langle\zeta,w\rangle|^p\,dS(\zeta)
=\frac{2\pi^N\,\Gamma(1+p/2)}{\Gamma(N+p/2)}\,|w|^p,\\
K_7&:=\int_{\mathbb B}|\langle z,w\rangle|^p\,(1-|z|^2)^q\,dV(z)
=\frac{\pi^N\,\Gamma(1+q)\,\Gamma(1+p/2)}{\Gamma(N+1+q+p/2)}\,|w|^p,\\
K_8&:=\int_{\mathbb B}|\langle z,w\rangle|^p\,dV(z)
=\frac{\pi^N\,\Gamma(1+p/2)}{\Gamma(N+1+p/2)}\,|w|^p.
\end{align*}
\begin{align*}
K_5'&:=\int_{\mathbb C^N}|\langle z,w\rangle|^p\,e^{-|z|^2}\,d\nu(z)
=N!\,(1)_{p/2}\,|w|^p,\\
K_6'&:=\int_{\mathbb S}|\langle\zeta,w\rangle|^p\,d\sigma(\zeta)
=\frac{(1)_{p/2}}{(N)_{p/2}}\,|w|^p,\\
K_7'&:=\int_{\mathbb B}|\langle z,w\rangle|^p\,(1-|z|^2)^q\,d\nu(z)
=\frac{(1)_q\,(1)_{p/2}}{(N+1)_{q+p/2}}\,|w|^p,\\
K_8'&:=\int_{\mathbb B}|\langle z,w\rangle|^p\,d\nu(z)
=\frac{(1)_{p/2}}{(N+1)_{p/2}}\,|w|^p.
\end{align*}
\begin{align*}
K_5''&:=\int_{\mathbb C^N}|\langle z,w\rangle|^{2m}\,e^{-|z|^2}\,dV(z)
=\pi^N\,m!\,|w|^{2m},\\
K_6''&:=\int_{\mathbb S}|\langle\zeta,w\rangle|^{2m}\,dS(\zeta)
=\frac{2\pi^N\,m!}{(N-1+m)!}\,|w|^{2m},\\
K_7''&:=\int_{\mathbb B}|\langle z,w\rangle|^{2m}\,(1-|z|^2)^k\,dV(z)
=\frac{\pi^N\,k!\,m!}{(N+k+m)!}\,|w|^{2m},\\
K_8''&:=\int_{\mathbb B}|\langle z,w\rangle|^{2m}\,dV(z)
=\frac{\pi^N\,m!}{(N+m)!}\,|w|^{2m}.
\end{align*}
\begin{align*}
K_5'''&:=\int_{\mathbb C^N}|\langle z,w\rangle|^{2m}\,e^{-|z|^2}\,d\nu(z)
=N!\,m!\,|w|^{2m},\\
K_6'''&:=\int_{\mathbb S}|\langle\zeta,w\rangle|^{2m}\,d\sigma(\zeta)
=\frac{(N-1)!\,m!}{(N-1+m)!}\,|w|^{2m},\\
K_7'''&:=\int_{\mathbb B}|\langle z,w\rangle|^{2m}\,(1-|z|^2)^k\,d\nu(z)
=\frac{N!\,k!\,m!}{(N+k+m)!}\,|w|^{2m},\\
K_8'''&:=\int_{\mathbb B}|\langle z,w\rangle|^{2m}\,d\nu(z)
=\frac{N!\,m!}{(N+m)!}\,|w|^{2m}.
\end{align*}

\section{Asymptotics}

Using \eqref{growth}, it is straightforward to compute the asymptotic growth
rates of the exact values of the integrals we have obtained above.
We list some below.
\medskip
\begin{equation*}
\left.
\begin{aligned}
J_3\sim J_3'&\sim q^{-(n+|\alpha|p)/2}\\
J_7\sim J_7'&\sim q^{-(n+p)/2}\,|y|^p
\end{aligned}
\right\}\qquad(q\to\infty),
\end{equation*}
\medskip
\begin{equation*}
\left.
\begin{aligned}
J_6\sim J_6'&\sim p^{-(n-1)/2}\,|y|^p\\
J_7\sim J_7'&\sim p^{-(n+1+2q)/2}\,|y|^p\\
J_8\sim J_8'&\sim p^{-(n+1)/2}\,|y|^p
\end{aligned}
\right\}\qquad(p\to\infty),
\end{equation*}
\medskip
\begin{equation*}
\left.
\begin{aligned}
K_3\sim K_3'&\sim q^{-(N+|\alpha|p/2)}\\
K_7\sim K_7'&\sim q^{-(N+p/2)}\,|w|^p
\end{aligned}
\right\}\qquad(q\to\infty),
\end{equation*}
\medskip
\begin{equation*}
\left.
\begin{aligned}
K_6\sim K_6'&\sim p^{-(N-1)}\,|w|^p\\
K_7\sim K_7'&\sim p^{-(N+q)}\,|w|^p\\
K_8\sim K_8'&\sim p^{-N}\,|w|^p
\end{aligned}
\right\}\qquad(p\to\infty).
\end{equation*}
\medskip

\bibliographystyle{amsalpha}

\end{document}